\documentclass{ifacconf}
\usepackage{natbib,amsmath,mathptmx}
\usepackage[pdftex]{graphicx}
  \graphicspath{{./pics/}}
  \DeclareGraphicsExtensions{.pdf}

\newcommand{\eps}{\varepsilon}
\newcommand{\supC}{\mbox{$\sup {\rm C}$}}
\newcommand{\supCC}{\mbox{$\sup {\rm cC}$}}

\begin{document}

\begin{frontmatter}
\title{On Algorithms and Extensions of Coordination Control of Discrete-Event Systems}

\author[CAS]{Jan~Komenda},
\author[CAS]{Tom{\' a}{\v s}~Masopust},
\author[CWI]{Jan~H.~van~Schuppen}

\address[CAS]{Institute of Mathematics, Academy of Sciences of the Czech Republic 
              {\v Z}i{\v z}kova 22, 616 62 Brno, Czech Republic\\ 
              (e-mails: komenda@ipm.cz, masopust@math.cas.cz)}
\address[CWI]{CWI, P.O. Box 94079, 1090 GB Amsterdam, The Netherlands\\ (e-mail: J.H.van.Schuppen@cwi.nl)}

\begin{abstract}
  In this paper, we further develop the coordination control scheme for discrete-event systems based on the Ramadge-Wonham framework. The notions of conditional decomposability, conditional controllability, and conditional closedness are revised and simplified, supremal conditionally controllable sublanguages of general non-prefix-closed languages are discussed, and a procedure for the computation of a coordinator for nonblockingness is presented.
\end{abstract}
\begin{keyword}
  Discrete-event system \sep supervisory control \sep coordination control \sep nonblockingness.
\end{keyword}
\end{frontmatter}

\section{Introduction}
  A distributed discrete-event system with synchronous communication is modeled as a parallel composition of two or more subsystems. Each subsystem has its own observation channel. The local control synthesis then consists in synthesizing local nonblocking supervisors for each subsystem. 
  
  Recently, \cite{KvS08} have proposed a coordination control architecture as a trade-off between the purely local control synthesis, which does not work in general, and the global control synthesis, which is not always possible because of complexity reasons. The coordination control approach has been developed for prefix-closed languages in~\cite{scl2011,automatica2011} and partially discussed for non-prefix-closed languages in~\cite{ifacwc2011}. A coordination control plug-in handling the case of prefix-closed languages has recently been implemented for libFAUDES, see~\cite{faudes}.
  
  In this paper, we further develop the coordination control scheme for discrete-event systems based on the Ramadge-Wonham framework. The notions of conditional decomposability, conditional controllability, and conditional closedness are revised and simplified, supremal conditionally controllable sublanguages of general non-prefix-closed languages are discussed, and a procedure for the computation of a coordinator for nonblockingness is presented.

  The paper is organized as follows. Section~\ref{sec:preliminaries} recalls the basic theory and revises the basic concepts. Section~\ref{sec:controlsynthesis} formulates the problem of coordination supervisory control. Section~\ref{sec:procedure} provides new results concerning non-prefix-closed languages, and Section~\ref{sec:nonblocking} discusses the construction of a nonblocking coordinator. Section~\ref{sec:revisited} revises the prefix-closed case, and Section~\ref{sec:conclusion} concludes the paper.

\section{Preliminaries and definitions}\label{sec:preliminaries}
  In this paper, we assume that the reader is familiar with supervisory control of discrete-event systems, where discrete-event systems are modeled as deterministic finite automata with partial transition functions, see~\cite{CL08}. 
  
  Let $E$ be a finite, nonempty set (of {\em events}), then $E^*$ denotes the set of all finite words over $E$; the {\em empty word\/} is denoted by $\eps$. A {\em generator\/} over $E$ is a construct $G=(Q,E,f,q_0,Q_m)$, where $Q$ is a finite set of {\em states}, $f: Q \times E \to Q$ is a {\em partial transition function}, $q_0 \in Q$ is the {\em initial state}, and $Q_m\subseteq Q$ is the set of {\em marked states}. In the usual way, $f$ can be extended to a function from $Q \times E^*$ to $Q$ by induction. The behavior of $G$ is described in terms of languages. The language {\em generated\/} by $G$ is the set $L(G) = \{s\in E^* \mid f(q_0,s)\in Q\}$, and the language {\em marked\/} by $G$ is the set $L_m(G) = \{s\in E^* \mid f(q_0,s)\in Q_m\}$.

  We restrict our attention to regular languages. A {\em (regular) language\/} $L$ over $E$ is a set $L\subseteq E^*$ such that there exists a generator $G$ with $L_m(G)=L$. The prefix closure of $L$ is the set $\overline{L}=\{w\in E^* \mid \exists u\in E^*, wu\in L\}$; $L$ is {\em prefix-closed\/} if $L=\overline{L}$.

  A {\em controlled generator\/} over $E$ is a structure $(G,E_c,\Gamma)$, where $G$ is a generator over $E$, $E_c \subseteq E$ is the set of {\em controllable events}, $E_{u} = E \setminus E_c$ is the set of {\em uncontrollable events}, and $\Gamma = \{\gamma \subseteq E \mid E_{u} \subseteq \gamma\}$ is a {\em set of control patterns}. A {\em supervisor\/} for the controlled generator $(G,E_c,\Gamma)$ is a map $S:L(G) \to \Gamma$. The {\em closed-loop system\/} associated with the controlled generator $(G,E_c,\Gamma)$ and the supervisor $S$ is defined as the minimal language $L(S/G)$ such that (i) $\eps \in L(S/G)$, and (ii) if $s \in L(S/G)$, $sa\in L(G)$, and $a \in S(s)$, then $sa \in L(S/G)$. We define $L_m(S/G) = L(S/G)\cap L_m(G)$. The supervisor disables transitions of $G$, but it cannot disable a transition with an uncontrollable event. If the closed-loop system is nonblocking, i.e., $\overline{L_m(S/G)}=L(S/G)$, then the supervisor $S$ is called {\em nonblocking}.

  Given a specification language $K$, the control objective of supervisory control is to find a nonblocking supervisor $S$ so that $L_m(S/G)=K$. For the monolithic case, such a supervisor exists if and only if $K$ is {\em controllable\/} with respect to $L(G)$ and $E_u$, that is, $\overline{K}E_u\cap L\subseteq \overline{K}$, and $K$ is {\em $L_m(G)$-closed}, that is, $K = \overline{K}\cap L_m(G)$. For uncontrollable specifications, controllable sublanguages are considered.  In this paper, $\supC(K,L,E_u)$ denotes the supremal controllable sublanguage of $K$ with respect to $L$ and $E_u$, which always exists and equals to the union of all controllable sublanguages of $K$, see~\cite{Won04}.

  A {\em projection\/} $P: E^* \to E_0^*$, $E_0\subseteq E$, is a homomorphism defined so that $P(a)=\eps$, for $a\in E\setminus E_0$, and $P(a)=a$, for $a\in E_0$. The {\em inverse image\/} of $P$ is denoted by $P^{-1}:E_0^* \to 2^{E^*}$. For $E_i$, $E_j$, $E_\ell \subseteq E$, we use the notation $P^{i+j}_{\ell}$ to denote the projection from $(E_i\cup E_j)^*$ to $E_\ell^*$. If $E_i\cup E_j=E$, we write only $P_\ell$. Moreover, $E_{i,u}=E_i\cap E_u$ denotes the sets of locally uncontrollable events.

  The synchronous product of languages $L_1\subseteq E_1^*$ and $L_2\subseteq E_2^*$ is defined by $L_1\| L_2=P_1^{-1}(L_1) \cap P_2^{-1}(L_2) \subseteq (E_1\cup E_2)^*$, where $P_i: (E_1\cup E_2)^*\to E_i^*$, $i=1,2$, are projections. For generators $G_1$ and $G_2$, the definition can be found in~\cite{CL08}. It holds that $L(G_1 \| G_2) = L(G_1) \| L(G_2)$ and $L_m(G_1 \| G_2)= L_m(G_1) \| L_m(G_2)$. In the automata framework, where the supervisor $S$ has a finite representation as a generator, the closed-loop system is a synchronous product of the supervisor and the plant. Thus, we can write $L(S/G)=L(S) \| L(G)$.

  Generators $G_1$ and $G_2$ are {\em conditionally independent\/} with respect to a generator $G_k$ if $E_r(G_1 \| G_2) \cap E_r(G_1) \cap E_r(G_2) \subseteq E_r(G_k)$, where for a generator $G$ over $E$, $E_r(G)=\{a\in E \mid \exists u,v\in E^*,\, uav\in L(G)\}$ is the set of all events appearing in words of $L(G)$. In other words, there is no simultaneous move in both $G_1$ and $G_2$ without the coordinator $G_k$ being also involved. From the practical viewpoint, we omit the element $E_r(G_1 \| G_2)$ because we do not want to compute the global plant $G_1\|G_2$.

  Now, the notion of decomposability is weakened. Moreover, it is simplified in comparison with our previous work, see~\cite{automatica2011}, but still equivalent. A language $K$ is {\em conditionally decomposable} with respect to event sets $E_1$, $E_2$, $E_k$ if $K = P_{1+k}(K)\| P_{2+k}(K)$. There always exists an extension of $E_k$ which satisfies the condition. The question which extension should be used (the minimal one?) requires further investigation. Polynomial-time algorithms for checking the condition and extending the event set are discussed in~\cite{CoRR12011733}.

  Languages $K$ and $L$ are {\em synchronously nonconflicting\/} if $\overline{K \| L} = \overline{K} \| \overline{L}$. Note that if $\overline{K}$ is conditionally decomposable, then the languages $P_{1+k}(K)$ and $P_{2+k}(K)$ are synchronously nonconflicting because $\overline{K} \subseteq \overline{P_{1+k}(K)\|P_{2+k}(K)} \subseteq \overline{P_{1+k}(K)} \| \overline{P_{2+k}(K)} = \overline{K}$. The following example shows that there is no relation between the conditional decomposability of $K$ and $\overline{K}$ in general.
  \begin{exmp}\label{ex1}
    Let $E_{1}=\{a_1,b_1,a,b\}$, $E_{2}=\{a_2,b_2,a,b\}$, $E_k=\{a,b\}$ be event sets, and let $K=\{a_1a_2a,a_2a_1a,b_1b_2b,b_2b_1b\}$. Then, $P_{1+k}(K) = \{a_1a,b_1b\}$, $P_{2+k}(K) = \{a_2a,b_2b\}$, and $K = P_{1+k}(K)\| P_{2+k}(K)$. Notice that $a_1b_2\in\overline{P_{1+k}(K)} \| \overline{P_{2+k}(K)}$, but $a_1b_2\notin\overline{K}$, which means that $\overline{K}$ is not conditionally decomposable.
    On the other hand, consider the language $L=\{\eps,ab,ba,abc,bac\}\subseteq \{a,b,c\}^*$ with $E_{1}=\{a,c\}$, $E_{2}=\{b,c\}$, $E_k=\{c\}$. Then, $\overline{L} = \overline{P_{1+k}(L)} \| \overline{P_{2+k}(L)} = P_{1+k}(L) \| P_{2+k}(L)$, and it is obvious that $L\neq \overline{L}$. \hfill$\triangleleft$
  \end{exmp}

\section{Coordination control synthesis}\label{sec:controlsynthesis}
  In this section, we formulate the coordination control problem and revise the necessary and sufficient conditions of \cite{ifacwc2011,scl2011,automatica2011} under which the problem is solvable.
  
  \begin{prob}\label{problem:controlsynthesis}
    Consider generators $G_1$, $G_2$ over $E_1$, $E_2$, respectively, and a coordinator $G_k$ over $E_k$. Let $K \subseteq L_m(G_1 \| G_2 \| G_k)$ be a specification. Assume that generators $G_1$ and $G_2$ are conditionally independent with respect to the coordinator $G_k$, and that the specification language $K$ and its prefix-closure $\overline{K}$ are conditionally decomposable with respect to $E_1$, $E_2$, $E_k$. The aim of the coordination control synthesis is to determine nonblocking supervisors $S_1$, $S_2$, $S_k$ for the respective generators such that $L_m(S_k/G_k)\subseteq P_k(K)$, $L_m(S_i/ [G_i \| (S_k/G_k) ])\subseteq P_{i+k}(K)$, for $i=1,2$, and the closed-loop system with the coordinator satisfies
    \begin{align*}
      L_m(S_1/ [G_1 \| (S_k/G_k) ]) ~ \| ~ L_m(S_2/ [G_2 \| (S_k/G_k) ]) & = K\,.
    \end{align*}
    $\hfill\diamond$
  \end{prob}

  Note that then $L(S_1/ [G_1 \| (S_k/G_k) ]) \| L(S_2/ [G_2 \| (S_k/G_k) ])=\overline{K}$ because $\overline{K} = \overline{L_m(S_1/ [G_1 \| (S_k/G_k) ]) \| L_m(S_2/ [G_2 \| (S_k/G_k) ])} \subseteq L(S_1/ [G_1 \| (S_k/G_k) ]) \| L(S_2/ [G_2 \| (S_k/G_k) ]) \subseteq \overline{K}$, and if such supervisors exist, their synchronous product is a nonblocking supervisor for the global plant, cf.~\cite{ifacwc2011}.
  
  One of the possible methods how to construct a suitable coordinator $G_k$ has been discussed in the literature, see \cite{ifacwc2011,scl2011,automatica2011}.
  \begin{alg}[Construction of a coordinator]\label{alg}
    Let $G_1$ and $G_2$ be two subsystems over $E_1$ and $E_2$, respectively, and let $K$ be a specification language. Construct the event set $E_k$ and the coordinator $G_k$ as follows:
    \begin{enumerate}
      \item Set $E_k = E_1\cap E_2$.
      \item Extend $E_k$ so that $K$ and $\overline{K}$ are conditional decomposable.
      \item Define $G_k = P_k(G_1) \parallel P_k(G_2)$.
    \end{enumerate}
  \end{alg}
  So far, the only known condition ensuring that the projected generator is smaller than the original one is the observer property. Therefore, we might need to add step (2b) to extend $E_k$ so that $P_k$ is also an $L(G_i)$-observer, for $i=1,2$, cf. Definition~\ref{def:observer}.

\subsection{Conditional controllability}
  Conditional controllability was introduced in \cite{KvS08} and later studied in \cite{ifacwc2011,scl2011,automatica2011}. In this paper, we revise and simplify this notion.

  \begin{defn}\label{def:conditionalcontrollability}
    A language $K\subseteq L(G_1\|G_2\|G_k)$ is {\em conditionally controllable\/} for generators $G_1$, $G_2$, $G_k$ and uncontrollable event sets $E_{1,u}$, $E_{2,u}$, $E_{k,u}$ if
    \begin{enumerate}
      \item\label{cc1} $P_k(K)$ is controllable wrt $L(G_k)$ and $E_{k,u}$,
      \item\label{cc2} $P_{1+k}(K)$ is controllable wrt $L(G_1) \parallel \overline{P_k(K)}$ and $E_{1+k,u}$,
      \item\label{cc3} $P_{2+k}(K)$ is controllable wrt $L(G_2) \parallel \overline{P_k(K)}$ and $E_{2+k,u}$.
    \end{enumerate}
    where $E_{i+k,u}=(E_i\cup E_k)\cap E_u$, $i=1,2$.
  \end{defn}

  The following result shows that every conditionally controllable and conditionally decomposable language is controllable. 
  \begin{prop}\label{prop3}
    Let $G_i$ be a generator over $E_i$, $i=1,2,k$, and let $G=G_1\|G_2\|G_k$. Let $K\subseteq L_m(G)$ be such that $\overline{K}$ is conditionally decomposable wrt $E_1$, $E_2$, $E_k$, and conditionally controllable for generators $G_1$, $G_2$, $G_k$ and uncontrollable event sets $E_{1,u}$, $E_{2,u}$, $E_{k,u}$. Then, $K$ is controllable with respect to $L(G)$ and $E_u$.
  \end{prop}
  \begin{pf}
    As $\overline{P_{1+k}(K)}$ is controllable wrt $L(G_1) \| \overline{P_k(K)}$ and $E_{1+k,u}$, and $\overline{P_{2+k}(K)}$ is controllable wrt $L(G_2) \| \overline{P_k(K)}$ and $E_{2+k,u}$, Lemma~\ref{feng} implies that $\overline{K} = \overline{P_{1+k}(K)} \| \overline{P_{2+k}(K)}$ is controllable wrt $L(G_1) \| \overline{P_k(K)} \| L(G_2) \| \overline{P_k(K)} = L(G) \| \overline{P_k(K)}$ and $E_u$, where the equality is by the commutativity of the synchronous product and the fact that $\overline{P_k(K)}\subseteq L(G_k)$. As $\overline{P_k(K)}$ is controllable wrt $L(G_k)$ and $E_{k,u}$, by Definition~\ref{def:conditionalcontrollability}, $L(G) \| \overline{P_k(K)}$ is controllable wrt $L(G) \| L(G_k) = L(G)$ by Lemma~\ref{feng}. By Lemma~\ref{lem_trans}, $\overline{K}$ is controllable wrt $L(G)$ and $E_u$. However, this means that $K$ is controllable wrt $L(G)$ and $E_u$, which was to be shown.
  \qed\end{pf}
  On the other hand, controllability does not imply conditional controllability.
  \begin{exmp}
    Let $L(G)=\overline{\{au\}}\|\overline{\{bu\}}=\overline{\{abu,bau\}}$. Then $K=\{a\}$ is controllable wrt $L(G)$ and $E_u=\{u\}$. Both $K$ and $\overline{K}$ are conditionally decomposable wrt event sets $\{a,u\}$, $\{b,u\}$, and $\{u\}$, and $P_k(K)=\{\varepsilon\}$ is not controllable wrt $\{u\}$ and $\{u\}$.
  \hfill$\triangleleft$\end{exmp}

  However, if the observer and local control consistency (LCC) properties are satisfied, this implication also holds. To prove this, we need the following two definitions, cf.~\cite{SB11,WW96}, respectively.

  \begin{defn}
    Let $L\subseteq E^*$ be a prefix-closed language, and let $E_0\subseteq E$. The projection $P_0:E^*\to E_0^*$ is {\em locally control consistent} (LCC) with respect to $s\in L$ if for all $\sigma_u\in E_0\cap E_u$ such that $P_0(s)\sigma_u\in P_0(L)$, it holds that either there does not exist any $u\in (E\setminus E_0)^*$ such that $su\sigma_u\in L$, or there exists $u\in (E_u\setminus E_0)^*$ such that $su\sigma_u\in L$. The projection $P_0$ is LCC with respect to a language $L$ if $P_0$ is LCC for all $s\in L$.
  \end{defn}

  \begin{defn}\label{def:observer}
    The projection $P_k:E^* \to E_k^*$, where $E_k\subseteq E$, is an {\em $L$-observer} for a language $L\subseteq E^*$ if, for all words $t\in P_k(L)$ and $s\in \overline{L}$, $P_k(s)$ is a prefix of $t$ implies that there exists $u\in E^*$ such that $su\in L$ and $P_k(su)=t$.
  \end{defn}

  \begin{prop}
    Let $L\subseteq E^*$ be a prefix-closed language, and let $K\subseteq L$ be a language such that $K$ is controllable with respect to $L$ and $E_u$. If $P_i$ is an $L$-observer, for $i\in\{k,1+k,2+k\}$, and LCC for $L$, then $K$ is conditionally controllable.
  \end{prop}
  \begin{pf}
    (1) Let $s\in \overline{P_k(K)}$, $a\in E_{k,u}$, and $sa\in P_k(L)$. Then, there exists $w\in\overline{K}$ such that $P_k(w)=s$. By the observer property, there exists $u\in (E\setminus E_k)^*$ such that $wua\in L$ and $P_k(wua)=sa$. By LCC, there exists $u'\in (E_u\setminus E_k)^*$ such that $wu'a\in L$, that is, $wu'a\in\overline{K}$ by the controllability. Hence $sa\in\overline{P_k(K)}$.
    (2) Let $s\in \overline{P_{1+k}(K)}$, $a\in E_{1+k,u}$, and $sa\in L(G_1)\|\overline{P_k(K)}$. Then, there exists $w\in\overline{K}$ such that $P_{1+k}(w)=s$. By the observer property, there exists $u\in (E\setminus E_{1+k})^*$ such that $wua\in L$ and $P_{1+k}(wua)=sa$. By LCC, there exists $u'\in (E_u\setminus E_{1+k})^*$ such that $wu'a\in L$, that is, $wu'a\in\overline{K}$ by controllability. Hence $sa\in\overline{P_{1+k}(K)}$.
  \qed\end{pf}
  
  For a generator $G$ with $n$ states, the time and space complexity of the verification whether $P$ is an $L(G)$-observer is $O(n^2)$, see~\cite{pcl08}. An algorithm extending the event set to satisfy the property runs in time $O(n^3)$ and linear space. The most significant consequence of the observer property is the following theorem.
  
  \begin{thm}[\cite{wong98}]
    If a projection $P$ is an $L(G)$-observer, for a generator $G$, then the minimal generator for the language $P(L(G))$ has no more states than $G$.
  \end{thm}

\subsection{Conditionally closed languages}
  Analogously to the notion of $L_m(G)$-closed languages, we define the notion of conditionally closed languages.

  \begin{defn}\label{def:conditionalclosed}
    A language $\emptyset\neq K \subseteq E^*$ is {\em conditionally closed\/} for generators $G_1$, $G_2$, $G_k$ if
    \begin{enumerate}
      \item\label{ccl1} $P_k(K)$ is $L_m(G_k)$-closed,
      \item\label{ccl2} $P_{1+k}(K)$ is $L_m(G_1) \| P_k(K)$-closed,
      \item\label{ccl3} $P_{2+k}(K)$ is $L_m(G_2) \| P_k(K)$-closed.
    \end{enumerate}
  \end{defn}

  If $K$ is conditionally closed and conditionally controllable, then there exists a nonblocking supervisor $S_k$ such that $L_m(S_k/G_k)=P_k(K)$, which follows from the basic theorem of supervisory control applied to $P_k(K)$ and $L(G_k)$, see~\cite{CL08}.
  
  As noted in \cite[page 164]{CL08}, if $K\subseteq L_m(G)$ is $L_m(G)$-closed, then so is the supremal controllable sublanguage of $K$. However, this does not imply that $P_k(K)$ is $L_m(G_k)$-closed, for $G=G_1\|G_2\|G_k$ such that $G_k$ makes $G_1$ and $G_2$ conditionally independent.
  \begin{exmp}\label{ex2}
    Let $E_{1} = \{a_1,a\}$, $E_{2}=\{a_2,a\}$, $E_k=\{a\}$, and $K=\{a_1a_2a,a_2a_1a\}$. Then, $P_{1+k}(K) = \{a_1a\}$, $P_{2+k}(K) = \{a_2a\}$, $P_{k}(K) = \{a\}$, and $K = P_{1+k}(K)\| P_{2+k}(K)$. Define the generators $G_1$, $G_2$, $G_k$ so that $L_m(G_1)=P_{1+k}(K)$, $L_m(G_2)=P_{2+k}(K)$, and $L_m(G_k)=\overline{P_{k}(K)}=\{\eps,a\}$. Then, $L_m(G)=K$ and $K$ is $L_m(G)$-closed. However, $P_k(K)\subset \overline{P_k(K)}$ is not $L_m(G_k)$-closed.
    \hfill$\triangleleft$
  \end{exmp}

\subsection{Coordination control synthesis}
  The following theorem is a simplified version of a result presented without proof in~\cite{ifacwc2011}.
  
  \begin{thm}\label{th:controlsynthesissafety}
    Consider the setting of Problem \ref{problem:controlsynthesis}. There exist nonblocking supervisors $S_1$, $S_2$, $S_k$ such that
    \begin{align}\label{eq:controlsynthesissafety}
      L_m(S_1/[G_1 \| (S_k/G_k)]) ~ & \| ~ L_m(S_2/[G_2 \| (S_k/G_k)]) =  K
    \end{align}
    if and only if the specification language $K$ is both conditionally controllable wrt generators $G_1$, $G_2$, $G_k$ and event sets $E_{1,u}$, $E_{2,u}$, $E_{k,u}$, and conditionally closed wrt $G_1$, $G_2$, $G_k$.
  \end{thm}
  \begin{pf}
    Let $K$ satisfy the assumptions, and let $G=G_1 \| G_2 \| G_k$. As $K\subseteq L_m(G)$, $P_k(K) \subseteq L_m(G_k)$. By the assumption, $P_k(K)$ is $L_m(G_k)$-closed and controllable wrt $L(G_k)$ and $E_{k,u}$. By~\cite{RW87}, there exists a nonblocking supervisor $S_k$ such that $L_m(S_k/G_k) = P_k(K)$. As $P_{1+k}(K) \subseteq L_m(G_1\|G_k)$ and $P_{1+k}(K) \subseteq (P_{k}^{1+k})^{-1} P_k(K)$, we have $P_{1+k}(K) \subseteq L_m(G_1) \| P_k(K)$. These relations and the assumption that the system is conditionally controllable and conditionally closed imply the existence of a nonblocking supervisor $S_1$ such that $L_m(S_1/ [ G_1 \| (S_k/G_k) ]) = P_{1+k}(K)$. A similar argument shows that there exists a nonblocking supervisor $S_2$ such that $L_m(S_2/ [ G_2 \| (S_k/G_k) ]) = P_{2+k}(K)$. As the languages $K$ and $\overline{K}$ are conditionally decomposable, $L_m(S_1/[G_1 \| (S_k/G_k)]) \parallel L_m(S_2/[G_2 \| (S_k/G_k)]) = P_{1+k}(K) \| P_{2+k}(K)  = K$.

    To prove the converse implication, $P_k$, $P_{1+k}$, $P_{2+k}$ are applied to (\ref{eq:controlsynthesissafety}), which can be rewritten as $K = L_m(S_1 \| G_1 \| S_2 \| G_2 \| S_k \| G_k)$. Thus, $P_k(K) = P_k\left(L_m(S_1 \| G_1 \| S_2 \| G_2 \| S_k \| G_k)\right) \subseteq L_m(S_k \| G_k) = L_m(S_k/G_k)$. On the other hand, $L_m(S_k/G_k) \subseteq P_k(K)$, cf. Problem~\ref{problem:controlsynthesis}. Hence, by the basic controllability theorem, $P_k(K)$ is controllable wrt $L(G_k)$ and $E_{k,u}$, and $L_m(G_k)$-closed. As $E_{1+k}\cap E_{2+k}=E_k$, the application of $P_{1+k}$ to (\ref{eq:controlsynthesissafety}) and Lemma~\ref{lemma:Wonham} give that $P_{1+k}(K) \subseteq L_m(S_1/ [ G_1 \| (S_k/G_k)]) \subseteq P_{1+k}(K)$. Taking $G_1 \| (S_k/G_k)$ as a new plant, we get that $P_{1+k}(K)$ is controllable wrt $L(G_1 \| (S_k/G_k))$ and $E_{1+k,u}$, and that it is $L_m(G_1 \| (S_k/G_k))$-closed. The case of $P_{2+k}$ is analogous.
  \qed\end{pf}

\section{Supremal conditionally controllable sublanguages}\label{sec:procedure}
  Let $\supCC(K, L, (E_{1,u}, E_{2,u}, E_{k,u}))$ denote the supremal conditionally controllable sublanguage of $K$ with respect to $L=L(G_1\| G_2\| G_k)$ and sets of uncontrollable events $E_{1,u}$, $E_{2,u}$, $E_{k,u}$. The supremal conditionally controllable sublanguage always exists, cf.~\cite{scl2011} for the case of prefix-closed languages.
  \begin{thm}\label{existence}
    The supremal conditionally controllable sublanguage of a given language $K$ always exists and is equal to the union of all conditionally controllable sublanguages of $K$.
  \end{thm}
  \begin{pf}
    Let $I$ be an index set, and let $K_i$, for $i\in I$, be conditionally controllable sublanguages of $K\subseteq L(G_1\|G_2\|G_k)$. 
    To prove that $P_k(\cup_{i\in I} K_i)$ is controllable wrt $L(G_k)$ and $E_{k,u}$, note that
    $P_k\left(\cup_{i\in I} \overline{K_i}\right)E_{k,u} \cap L(G_k) 
      = \cup_{i\in I} \left(P_k(\overline{K_i})E_{k,u} \cap L(G_k) \right) 
      \subseteq \cup_{i\in I} P_k(\overline{K_i}) 
      = P_k\left(\cup_{i\in I}\overline{K_i}\right)$,
    where the inclusion is by controllability of $P_k(K_i)$ wrt $L(G_k)$ and $E_{k,u}$.
    Next, to prove that
    $
      P_{1+k}\left(\cup_{i\in I} \overline{K_i}\right)E_{1+k,u} 
      \cap L(G_1)\|P_k\left(\cup_{i\in I} \overline{K_i}\right)
      \subseteq P_{1+k}\left(\cup_{i\in I} \overline{K_i}\right),
    $
    note that $P_{1+k}\left(\cup_{i\in I} \overline{K_i}\right)E_{1+k,u} \cap L(G_1)\|P_k\left(\cup_{i\in I} \overline{K_i}\right)$
    \begin{align*}
      &= \cup_{i\in I} \left( P_{1+k}(\overline{K_i})E_{1+k,u}\right) 
        \cap \cup_{i\in I} \left(L(G_1)\|P_k(\overline{K_i})\right)\\
      &= \cup_{i\in I}\cup_{j\in I} \left( P_{1+k}(\overline{K_i})E_{1+k,u} 
        \cap L(G_1)\|P_k(\overline{K_j})\right)\,.
    \end{align*}
    Consider different indexes $i,j\in I$ such that
    $
      P_{1+k}(\overline{K_i})E_{1+k,u} 
      \cap L(G_1)\|P_k(\overline{K_j})
      \not\subseteq P_{1+k}\left(\cup_{i\in I} \overline{K_i}\right).
    $
    Then, there exist $x\in P_{1+k}(\overline{K_i})$ and $u\in E_{1+k,u}$ such that $xu\in L(G_1)\|P_k(\overline{K_j})$, and $xu \notin P_{1+k}\left(\cup_{i\in I} \overline{K_i}\right)$. It follows that
      $P_k(x)\in P_k(\overline{K_i})$ and
      $P_k(xu)\in P_k(\overline{K_j})$.
    If $P_k(xu)\in P_k(\overline{K_i})$, then $xu\in L(G_1)\|P_k(\overline{K_i})$, and controllability of $P_{1+k}(\overline{K_i})$ wrt $L(G_1)\|P_k(\overline{K_i})$ implies that $xu\in P_{1+k}\left(\cup_{i\in I} \overline{K_i}\right)$; hence $P_k(xu)\notin P_k(\overline{K_i})$. If $u\notin E_{k,u}$, then $P_k(xu)=P_k(x)\in P_k(\overline{K_i})$, which is not the case. Thus, $u\in E_{k,u}$. As $P_k(\overline{K_i})\cup P_k(\overline{K_j})\subseteq L(G_k)$, we get that $ P_k(xu) = P_k(x)u \in L(G_k)$. However, controllability of $P_k(\overline{K_i})$ wrt $L(G_k)$ and $E_{k,u}$ implies that $ P_k(xu) \in P_k(\overline{K_i})$. This is a contradiction.
    As the case for $P_{2+k}$ is analogous, the proof is complete.
  \qed\end{pf}  
  
  Consider the setting of Problem~\ref{problem:controlsynthesis}, and define the languages
  \begin{equation}\tag{*}\label{eq0}
    \boxed{
    \begin{aligned}
      \supC_k     & =  \supC(P_k(K), L(G_k), E_{k,u})\,,\\
      \supC_{1+k} & =  \supC(P_{1+k}(K), L(G_1) \| \overline{\supC_k}, E_{1+k,u})\,,\\
      \supC_{2+k} & =  \supC(P_{2+k}(K), L(G_2) \| \overline{\supC_k}, E_{2+k,u})\,.
    \end{aligned}}
  \end{equation}

  The following inclusion always holds.
  \begin{lem}\label{lemma16}
    Consider the setting of Problem~\ref{problem:controlsynthesis}, and the languages defined in (\ref{eq0}). Then, $P_k(\supC_{i+k}) \subseteq \supC_k$, for $i=1,2$.
  \end{lem}
  \begin{pf}
    By definition, $P_k(\supC_{i+k}) \subseteq \overline{\supC_k}$ and $P_k(\supC_{i+k}) \subseteq P_k(K)$.
    To prove that $\overline{\supC_k}\cap P_k(K)$ is a subset of $\supC_k$, it is sufficient to show that $\overline{\supC_k}\cap P_k(K)$ is controllable with respect to $L(G_k)$ and $E_{k,u}$.
    Thus, assume that $s\in \overline{\overline{\supC_k}\cap P_k(K)}$, $u\in E_{k,u}$, and $su\in L(G_k)$. 
    By controllability of $\supC_k$, $su\in\overline{\supC_k}\subseteq \overline{P_k(K)}$, that is, there exists $v$ such that $suv\in\supC_k\subseteq P_k(K)$. This means that $suv\in \overline{\supC_k}\cap P_k(K)$, which implies that $su\in \overline{\overline{\supC_k}\cap P_k(K)}$. This completes the proof.
  \qed\end{pf}

  If also the opposite inclusion holds, then we immediately have the supremal conditionally-controllable sublanguage.
  \begin{thm}\label{thm2}
    Consider the setting of Problem~\ref{problem:controlsynthesis}, and the languages defined in (\ref{eq0}). If $\supC_k\subseteq P_k(\supC_{i+k})$, for $i=1,2$, then $\supC_{1+k} \| \supC_{2+k} = \supCC(K, L, (E_{1,u}, E_{2,u}, E_{k,u}))$.
  \end{thm}
  \begin{pf}
    Let $\supCC = \supCC(K, L, (E_{1,u}, E_{2,u}, E_{k,u}))$ and $M = \supC_{1+k} \| \supC_{2+k}$.
    To prove $M\subseteq \supCC$, we show that (i) $M \subseteq K$ and (ii) $M$ is conditionally controllable wrt $G_1$, $G_2$, $G_k$ and $E_{1,u}$, $E_{2,u}$, $E_{k,u}$. To this aim, $M = \supC_{1+k} \| \supC_{2+k} \subseteq P_{1+k}(K) \| P_{2+k}(K) = K$ because $K$ is conditionally decomposable. Moreover, $P_k(M) = P_k(\supC_{1+k}) \cap P_k(\supC_{2+k})=\supC_k$, which is controllable wrt $L(G_k)$ and $E_{k,u}$. Similarly, $P_{i+k}(M) = \supC_{i+k} \| P_k(\supC_{j+k}) = \supC_{i+k} \| \supC_{k} = \supC_{i+k}$, for $j\neq i$, by Lemma~\ref{lemma16}, which is controllable wrt $L(G_i)\|\overline{P_k(M)}$. Hence, $M\subseteq \supCC$.

    To prove the opposite inclusion, by Lemma~\ref{lem11}, it is sufficient to show that $P_{i+k}(\supCC)\subseteq \supC_{i+k}$, for $i=1,2$. To prove this $P_{1+k}(\supCC)$ is controllable wrt $L(G_1)\| \overline{P_k(\supCC)}$ and $E_{1+k,u}$, and $L(G_1) \| \overline{P_k(\supCC)}$ is controllable wrt $L(G_1) \| \overline{\supC_k}$ and $E_{1+k,u}$ by Lemma~\ref{feng} because $P_k(\supCC)$ being controllable wrt $L(G_k)$ implies it is controllable wrt $\overline{\supC_k} \subseteq L(G_k)$ and $E_{k,u}$. By Lemma~\ref{lem_trans}, $P_{1+k}(\supCC)$ is controllable wrt $L(G_1) \| \overline{\supC_k}$ and $E_{1+k,u}$, which implies that $P_{1+k}(\supCC)\subseteq \supC_{1+k}$. The other case is analogous. Hence, $\supCC \subseteq M$ and the proof is complete.
  \qed\end{pf}

  \begin{exmp}
    This example shows that the inclusion $\supC_k\subseteq P_k(\supC_{i+k})$ does not hold in general. Moreover, it shows that it does not hold even if the projections are observers or satisfy the LCC property. Consider two systems $G_1$, $G_2$, and the specification $K$ as shown in Fig.~\ref{figEx}.
    \begin{figure}[ht]
      \centering
      \includegraphics[scale=.7]{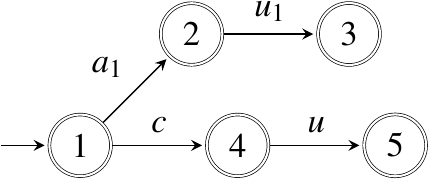}
      \qquad
      \includegraphics[scale=.7]{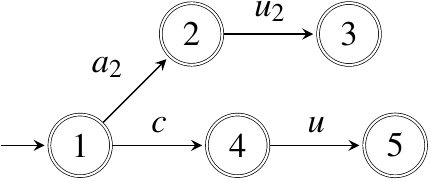}
      \qquad
      \includegraphics[scale=.7]{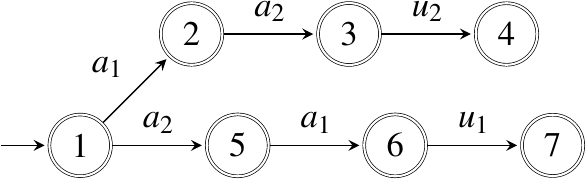}
      \caption{Generators $G_1$, $G_2$, and the specification.}
      \label{figEx}
    \end{figure}
    The controllable events are $E_c=\{a_1,a_2,c\}$, and the coordinator events are $E_k=\{a_1,a_2,c,u\}$. Construct the coordinator $G_k=P_k(G_1)\| P_k(G_2)$. It can be verified that $K$ is conditionally decomposable, $\supC_k=\overline{\{a_1a_2,a_2a_1\}}$, $\supC_{1+k}=\overline{\{a_2a_1u_1\}}$, and $\supC_{2+k}=\overline{\{a_1a_2u_2\}}$. Hence, $\supC_k\not\subseteq P_k(\supC_{i+k})$.
    It can also be verified that the projections $P_k$, $P_{1+k}$, $P_{2+k}$ are $L(G_1\|G_2)$-observers and LCC for $L(G_1\|G_2)$.
  \hfill$\triangleleft$\end{exmp}

  \begin{prop}
    Consider the languages of (\ref{eq0}). Let the number of states of the supervisor $\supC_k$ be $n$ and the number of states of supervisors $\supC_{i+k}$ be $n_i$. There is an $O(n\cdot n_i)$ algorithm deciding whether $\supC_k\subseteq P_k(\supC_{i+k})$, for $i=1,2$.
  \end{prop}
  \begin{pf}
    Consider a nondeterministic finite automaton, cf.~\cite{Sipser}, for the language $P_k(\supC_{i+k})$ constructed from the generator for $\supC_{i+k}$ by replacing projected events with $\eps$, and a deterministic finite automaton for the complement of $\supC_k$. These automata are constructed in time linear wrt the number of states. To verify that $P_k(\supC_{i+k}) \cap \mbox{co-}(\supC_k) = \emptyset$ by checking reachability of a marked state in the product automaton takes time $O(n\cdot n_i)$; here ``co-'' stands for the complement.
  \qed\end{pf}

  Note that if we have any specification $K$ which is conditionally decomposable, then the specification $K\|L$ is also conditionally decomposable. The opposite is not true.
  \begin{lem}
    Let $K$ be conditionally decomposable with respect to event sets $E_1$, $E_2$, $E_k$, and let $L=L_1\|L_2\|L_k$, where $L_i\subseteq E_i^*$, for $i=1,2,k$. Then, $K\|L$ is conditionally decomposable with respect to event sets $E_1$, $E_2$, $E_k$.
  \end{lem}

  \begin{exmp}\label{ex3}
    Database transactions are examples of discrete-event systems that need to be controlled to avoid incorrect behaviors. Our model of a transaction to the database is a sequence of request ($r$), access ($a$), and exit ($e$) operations. Usually, several (but a limited number of) users access the database, which can lead to inconsistencies when executed concurrently because not all the interleavings of operations give a correct behavior. We consider the case of three users with events $r_i,a_i,e_i$, $i=1,2,3$. All possible schedules are given by the language of the plant $G=G_1\|G_2\|G_3$ over the event set $E=\{r_1,r_2,r_3,a_1,a_2,a_3,e_1,e_2,e_3\}$, where $G_1$, $G_2$, $G_3$ are defined as in Fig.~\ref{figA}, and the set of controllable events is $E_c=\{a_1,a_2,a_3\}$.%
    \begin{figure}
      \centering
      \includegraphics[scale=.7]{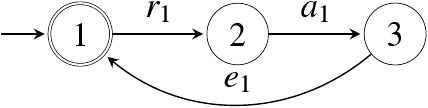}
      \qquad
      \includegraphics[scale=.7]{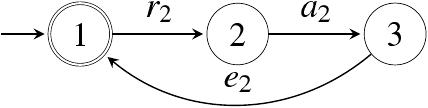}
      \qquad
      \includegraphics[scale=.7]{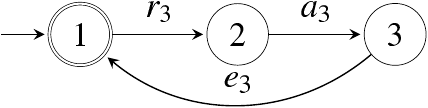}
      \caption{Generators $G_i$, $i=1,2,3$.}
      \label{figA}
    \end{figure}
    \begin{figure}
      \centering
      \includegraphics[scale=.7]{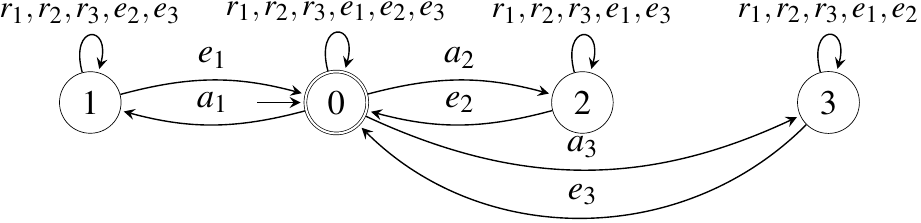}
      \caption{The specification $K$.}
      \label{figspecK}
    \end{figure}
    The specification language $K$, depicted in Fig.~\ref{figspecK}, describes the correct behavior consisting in finishing the transaction in the exit stage before another transaction can proceed to the exit phase.
    \begin{figure}
      \centering
      \includegraphics[scale=.7]{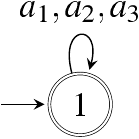}
      \caption{The coordinator $G_k$, where $supC_{k}=G_k$.}
      \label{figC}
    \end{figure}
    For $E_k=\{a_1,a_2,a_3\}$ and the coordinator $G_k = P_{k}(G_1) \| P_{k}(G_2) \| P_k(G_3)$, we can compute $\supC_k$, see Fig.~\ref{figC},  and $\supC_{1+k}$, $\supC_{2+k}$, $\supC_{3+k}$, Fig.~\ref{supsfig}, and to verify that the assumptions of Theorem~\ref{thm2} are satisfied. 
    \begin{figure}
      \centering
      \includegraphics[scale=.7]{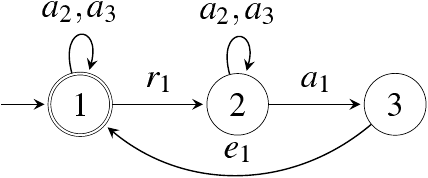}\qquad
      \includegraphics[scale=.7]{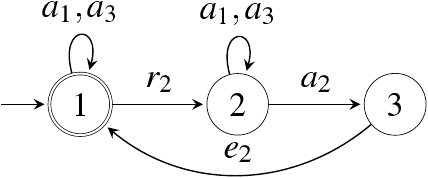}\qquad
      \includegraphics[scale=.7]{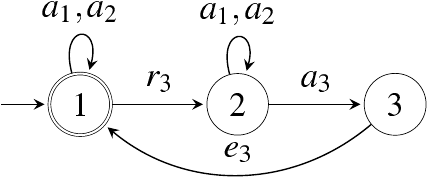}
      \caption{Supervisors $\supC_{1+k}$, $\supC_{2+k}$, and $\supC_{3+k}$.}
      \label{supsfig}
    \end{figure}  
    The solution is optimal: the supremal conditionally-controllable sublanguage of $K$ coincides with the supremal controllable sublanguage of $K$. Moreover, independently on the size of the global plant, the local supervisors have only three states.
  \hfill$\triangleleft$\end{exmp}

\section{Coordinator for nonblockingness}\label{sec:nonblocking}
   So far, we have only considered the coordinator for safety. In this section, we discuss the coordinator for nonblockingness. To this end, we first prove a fundamental theoretical result and then give an algorithm for the construction of a coordinator for nonblockingness. 
   
   Recall that a generator $G$ is nonblocking if $\overline{L_m(G)}=L(G)$.
  \begin{thm}\label{thm22}
    Consider languages $L_1\subseteq E_1^*$ and $L_2\subseteq E_2^*$, and let $P_0:(E_1\cup E_2)^*\to E_0^*$, with $E_1\cap E_2\subseteq E_0$, be an $L_i$-observer, for $i=1,2$. Let $G_0$ be a nonblocking generator with $L_m(G_0)=P_0(L_1)\|P_0(L_2)$. Then $\overline{L_1\|L_2\|L_m(G_0)} = \overline{L_1}\|\overline{L_2}\|\overline{L_m(G_0)}$, that is, the system is nonblocking.
  \end{thm}
  \begin{pf}
    Let $\overline{L} = \overline{L_1\|L_2\|L_0} = \overline{(L_1\|L_0)\|(L_2\|L_0)}$. By Lemma~\ref{fengT41}, $\overline{(L_1\|L_0)\|(L_2\|L_0)}=\overline{(L_1\|L_0)}\|\overline{(L_2\|L_0)}$ if and only if it holds $\overline{P_0(L_1\|L_0)\|P_0(L_2\|L_0)}=\overline{P_0(L_1\|L_0)}\|\overline{P_0(L_2\|L_0)}$, because if $P_0$ is an $L_i$-observer, $i=1,2$, and $P_0$ is an $L_0$-observer, $P_0$ is also an $L_i\|L_0$-observer by \cite{pcl06}. However, for our choice of the coordinator, this equality always holds because $\overline{P_0(L_1\|L_0)\|P_0(L_2\|L_0)} = \overline{L_0}$, and $\overline{P_0(L_1\|L_0)}\|\overline{P_0(L_2\|L_0)}=\overline{L_0}\|\overline{L_0}=\overline{L_0}$. It remains to show that $\overline{L_i\|L_0}=\overline{L_i}\|\overline{L_0}$, for $i=1,2$. Using Lemma~\ref{fengT41} again, we get that this holds if and only if $\overline{P_0(L_i\|L_0)}=\overline{P_0(L_i)}\|\overline{L_0}$. This always holds because $\overline{P_0(L_i\|L_0)}=\overline{L_0}$, and $\overline{P_0(L_i)}\|\overline{L_0} = \overline{P_0(L_i)}\|\overline{P_0(L_1)\|P_0(L_2)}=\overline{P_0(L_1)\|P_0(L_2)}=\overline{L_0}$ because $\overline{P_0(L_1)\|P_0(L_2)}\subseteq \overline{P_0(L_i)}$.
  \qed\end{pf}
  
  Hence, for supervisors $\supC_{1+k}$ and $\supC_{2+k}$, we choose \[C=P_0(\supC_{1+k}) \| P_0(\supC_{2+k})\,,\] for the projection $P_0$ being a $\supC_{i+k}$-observer, for $i=1,2$. Then, by Theorem~\ref{thm22},
  \begin{align*}
    \overline{\supC_{1+k}\|\supC_{2+k}\|C} &= \overline{\supC_{1+k}\|\supC_{2+k}}\\ 
    &= \overline{\supC_{1+k}} \| \overline{\supC_{2+k}} \| \overline{C}\,,
  \end{align*}
  thus $C$ is the language of a non-blocking coordinator.
  
  \begin{alg}[Computation of a nonblocking coordinator] $ $\\
    Consider the notation above.
    \begin{enumerate}
      \item Compute $\supC_{1+k}$ and $\supC_{2+k}$ as defined in (\ref{eq0}).
      \item If the projection $P_k$ is not a $\supC_{1+k}$-observer or not a $\supC_{2+k}$-observer, extend the event set $E_k$ so that $P_k$ is both  a $\supC_{1+k}$- and a $\supC_{2+k}$-observer.
      \item Define the nonblocking coordinator as the minimal nonblocking generator for $C=P_k(\supC_{1+k})\|P_k(\supC_{2+k})$.
    \end{enumerate}
  \end{alg}

\section{Supremal prefix-closed languages}\label{sec:revisited}
  In this section, we revise the case of prefix-closed languages. Moreover, we use LCC instead of output control consistency (OCC), cf.~\cite{automatica2011}.
  
  \begin{thm}\label{thm2pc}
    Let $K\subseteq L=L(G_1\|G_2\|G_k)$ be a prefix-closed language, where $G_i$ is over $E_i$, $i=1,2,k$. Assume that $K$ is conditionally decomposable, and define $\supC_k$, $\supC_{1+k}$, $\supC_{2+k}$ as in (\ref{eq0}). Let $P^{i+k}_k$ be an $(P^{i+k}_i)^{-1}(L(G_i))$-observer and LCC for $(P^{i+k}_i)^{-1}(L(G_i))$, $i=1,2$. Then, $\supC_{1+k} \| \supC_{2+k} = \supCC(K, L, (E_{1,u}, E_{2,u}, E_{k,u}))$.
  \end{thm}
  \begin{pf}
    Denote $\supCC = \supCC(K, L, (E_{1,u}, E_{2,u}, E_{k,u}))$, $M =\supC_{1+k} \| \supC_{2+k}$. It is shown in~\cite{automatica2011} that $\supCC \subseteq M$ and $M \subseteq K$. To prove $P_k(M) E_{k,u} \cap L(G_k) \subseteq P_k(M)$, let $x\in P_k(M)$ and $a\in E_{k,u}$ be such that $xa\in L(G_k)$. To show $xa \in P_k(M)= P^{1+k}_k(\supC_{1+k}) \cap P^{2+k}_k(\supC_{2+k})$, there exists $w\in M$ such that $P_k(w)=x$, and it is shown in~\cite{automatica2011} that there exists $u\in (E_1\setminus E_k)^*$ such that $P_{1+k}(w)ua\in (P^{1+k}_1)^{-1}(L(G_1))$ and $P_{1+k}(w)\in L(G_1) \| \supC_k$. As $P^{1+k}_k$ is LCC for $(P^{1+k}_1)^{-1}(L(G_1))$, there exists $u'\in (E_u\setminus E_k)^*$ such that $P_{1+k}(w)u'a\in (P^{1+k}_1)^{-1}(L(G_1))$. The controllability of $\supC_{1+k}$ then implies $P_{1+k}(w)u'a \in \supC_{1+k}$, i.e., $xa \in P^{1+k}_k(\supC_{1+k})$. Analogously, $xa\in P^{2+k}_k(\supC_{2+k})$. Thus, $xa\in P_k(M)$. The rest of the proof is the same as in~\cite{automatica2011}.
  \qed\end{pf}

  The conditions of Theorem~\ref{thm2pc} imply that $P_k$ is LCC for $L$.
  \begin{lem}\label{lem46}
    Let $L(G_i)\subseteq E_i^*$, $i=1,2$, $E=E_1\cup E_2$, and let $P_i: E^*\to E_i^*$, $i=1,2,k$ and $E_k\subseteq E$, be projections. If $E_1\cap E_2 \subseteq E_k$ and $P^{i+k}_k$ is LCC for $(P^{i+k}_i)^{-1}(L(G_i))$, $i=1,2$, then $P_k$ is LCC for $L=L(G_1\|G_2\|G_k)$.
  \end{lem}
  \begin{pf}
    For $s\in L$ and $\sigma_u\in E_{k,u}$, assume that there exists $u\in (E\setminus E_k)^*$ such that $su\sigma_u \in L$. Then, $P_{i+k}(su\sigma_u) = P_{i+k}(s)P_{i+k}(u)\sigma_u\in (P^{i+k}_i)^{-1}(L(G_i))$ implies that there exists $v_{i}\in (E_{i+k,u}\setminus E_k)^*$, $i=1,2$, such that $P_{i+k}(s)v_{i}\sigma_u\in (P^{i+k}_i)^{-1}(L(G_i))$. As $P_k(v_{i})=\eps$, $P_i(v_{i})=v_i$, we get $P_i(s)P_i(v_i)P_i(\sigma_u)\in L(G_i)$, $i=1,2,k$. Consider $u'\in \{v_1\}\|\{v_2\}$. Then $P_i(u')=v_i$ and, thus, $su'\sigma_u \in L$. Moreover, $u'\in (E_u\setminus E_k)^*$.
  \qed\end{pf}

  It is an open problem how to verify that $P_{i+k}$ is LCC for $L$ without computing the whole plant.
  \begin{thm}\label{thm4}
    Consider the setting of Theorem \ref{thm2pc}. If, in addition, $L(G_k)\subseteq P_k(L)$ and $P_{i+k}$ is LCC for $L$, for $i=1,2$, then $\supC(K, L, E_{u}) = \supCC(K, L, (E_{1,u}, E_{2,u}, E_{k,u}))$.
  \end{thm}
  \begin{pf}
    It was shown in~\cite{automatica2011} that $P_k$ is an $L$-observer. By Lemma~\ref{lem46}, $P_k$ is LCC for $L$. Denote $\supC = \supC(K,L,E_u)$. We prove that $P_k(\supC)$ is controllable wrt $L(G_k)$. Assume $t\in P_k(\supC)$, $a\in E_{k,u}$, and $ta\in L(G_k)\subseteq P_k(L)$. We proved in~\cite{automatica2011} that there exists $s\in\supC$ and $u\in (E\setminus E_k)^*$ such that $sua\in L$ and $P_k(sua)=ta$. By the LCC property of $P_k$, there exists $u'\in (E_{u}\setminus E_k)^*$ such that $su'a\in L$. By controllability of $\supC$ wrt $L$, $su'a\in\supC$, i.e., $P_k(su'a)=ta\in P_k(\supC)$. Thus, (1) of Definition~\ref{def:conditionalcontrollability} holds. By~\cite{automatica2011}, $P_{i+k}$ is an $L$-observer, for $i=1,2$. To prove (2) of Definition~\ref{def:conditionalcontrollability}, assume that $t\in P_{i+k}(\supC)$, $1\le i\le 2$, $a\in E_{i+k,u}$, and $ta\in L(G_i)\|P_k(\supC)$. We proved in~\cite{automatica2011} that there exists $s\in\supC$ and $u\in (E\setminus E_k)^*$ such that $sua\in L$ and $P_{i+k}(sua)=ta$. As $P_{i+k}$ is LCC for $L$, there exists $u'\in (E_{u}\setminus E_{1+k})^*$ such that $su'a\in L$. Then, the controllability of $\supC$ wrt $L$ implies that $su'a\in\supC$, that is, $P_{i+k}(su'a)=ta\in P_{i+k}(\supC)$. The other inclusion is the same as in~\cite{automatica2011}.
  \qed\end{pf}

\section{Conclusion}\label{sec:conclusion}
  We have revised, simplified, and extended the coordination control scheme for discrete-event systems. These results have been used, for the case of prefix-closed languages, in the implementation of the coordination control plug-in for libFAUDES. Note that a general procedure for the computation of supremal conditionally-controllable sublanguages is still missing. This requires further investigation.

\section*{Auxiliary results}\label{appendix}
	\begin{lem}[Proposition~4.6, \cite{FLT}]\label{feng}
		Let $L_i\subseteq E_i^*$, $i=1,2$, be prefix-closed languages, and let $K_i\subseteq L_i$ be controllable with respect to $L_i$ and $E_{i,u}$, $E=E_1\cup E_2$. If $K_1$ and $K_2$ are synchronously nonconflicting, then $K_1\|K_2$ is controllable with respect to $L_1\|L_2$ and $E_u$.
	\end{lem}

	\begin{lem}[\cite{automatica2011}]\label{lem_trans}
		Let $K\subseteq L \subseteq M$ be languages over $E$ such that $K$ is controllable with respect to $\overline{L}$ and $E_u$, and $L$ is controllable with respect to $\overline{M}$ and $E_u$. Then, $K$ is controllable with respect to $\overline{M}$ and $E_u$.
	\end{lem}

  \begin{lem}[\cite{Won04}]\label{lemma:Wonham}
    Let $P_k : E^*\to E_k^*$, $L_i\subseteq E_i^*$, $E_i\subseteq E$, $i=1,2$, $E_k\supseteq E_1\cap E_2$. Then, $P_k(L_1\| L_2)=P_k(L_1) \| P_k(L_2)$.
  \end{lem}

  \begin{lem}[\cite{automatica2011}]\label{lem11}
    Let $L_i\subseteq E_i^*$, $i=1,2$, and $P_i : (E_1\cup E_2)^* \to E_i^*$. Let $A\subseteq (E_1\cup E_2)^*$ be a language such that $P_1(A)\subseteq L_1$ and $P_2(A)\subseteq L_2$. Then $A\subseteq L_1\|L_2$.
  \end{lem}

  \begin{lem}[\cite{pcl06}]\label{fengT41}
    Let $L_i\subseteq E_i^*$, $i=1,2$, and let $E_1\cap E_2\subseteq E_0$. If $P_{i,0}:E_i^* \to (E_i\cap E_0)^*$ is an $L_i$-observer, $i=1,2$, then $\overline{L_1\|L_2} = \overline{L_1}\|\overline{L_2}$ iff $\overline{P_{1,0}(L_1)\|P_{2,0}(L_2)} = \overline{P_{1,0}(L_1)}\|\overline{P_{2,0}(L_2)}$.
  \end{lem}

\begin{ack}
  The research has been supported by the GA{\v C}R grants no. P103/11/0517 and P202/11/P028, and by RVO: 67985840.
\end{ack}

\bibliography{wodes2012}

\end{document}